\documentstyle[12pt]{article}
\titlepage

\title{Ricci Curvature and Gauss Maps of Minimal Submanifolds}
\author{Richard Atkins } 
		
\date{Subject classification: 53C42 
       Ricci curvature, Grassmannian}
\newtheorem{fact}{Fact}

\newtheorem{lemma}[fact]{Lemma}

\newtheorem{theorem}[fact]{Theorem}
\newtheorem{corollary}[fact]{Corollary}
\begin{document}
\maketitle
\begin{abstract} 
We present conditions on the Ricci curvature for complete, oriented, minimal
submanifolds of Euclidean space, as well as the standard unit sphere, when the Gauss maps
are bounded embeddings. 
\end{abstract}

\newpage

\section{Introduction}
A fundamental approach to the study of minimal submanifolds was inspired by Bernstein,
whose theorem relates properties about the Gauss map of a minimal surface to
information about the minimal surface itself. Herein we investigate upper bounds
on the Ricci curvature of minimal submanifolds of Euclidean space and the standard 
unit sphere
when the Gauss map is a bounded embedding. A {\it bounded embedding} is an embedding
whose image ${\cal I}$ in the oriented Grassmann manifold $\widetilde{G}$ is bounded 
with respect to the metric on ${\cal I}$ induced from the canonical metric on 
$\widetilde{G}$. We prove the following two theorems.
 \vspace*{0.1in} \\
{\bf Theorem A.} {\it Let $M$ be a complete, oriented, minimal submanifold of $R^{k}$ and 
suppose the 
Gauss map is a bounded embedding. Then $Ric \hspace{0.02in} < \hspace{0.02in}  0$ and
\begin{eqnarray}
\sup Ric & = & 0 \nonumber
\end{eqnarray} }
\\
{\bf Theorem B.} {\it Let $M^{m}$ be a complete, oriented, minimal submanifold of 
a unit sphere 
and suppose the second Gauss map is a bounded embedding. 
Then $Ric \hspace*{0.02in} < \hspace*{0.02in} m-1$. 
Furthermore, $M$ is compact  if and only if 
\begin{eqnarray}
\sup Ric & < & m-1 \nonumber
\end{eqnarray} }
\vspace*{0.1in} \\
Bounds for the Ricci curvature of complete, oriented, minimal hypersurfaces of the
unit sphere have been obtained by Hasanis and Vlachos \cite{dd}. A pinching theorem
for minimal submanifolds of a sphere having positive Ricci curvature was proved by 
Ejiri \cite{aa}.

\section{The First Gauss Map}

The oriented Grassmannian $\widetilde{G}_{m,n}$ is the Riemannian manifold of oriented 
$m$-planes in $R^{m+n}$ with canonical metric
$g_{c}$ determined as follows. Choose any two points $P,Q \in \widetilde{G}_{m,n}$ and let
$\xi_{1},...,\xi_{m}$ and $\zeta_{1},...,\zeta_{m}$ be oriented, orthonormal bases of $P$ 
and $Q$, respectively. Form the $m\times m$ matrix $\alpha$ by $\alpha_{ij} = 
<\xi_{i}, \zeta_{j}>$ and let $A$ be the product of $\alpha$ with
its transpose: $A:= \alpha \alpha^{T}$. $A$ is non-negative and symmetric with eigenvalues
$\lambda_{1}^{2},...,\lambda_{m}^{2}$, say, where  $0\leq \lambda_{i} \leq 1$. 
Define $d_{c}: \widetilde{G}_{m,n}\times\widetilde{G}_{m,n} \longrightarrow 
[0,\sqrt{m}\pi/2]$ by
\begin{eqnarray}
d_{c}(P,Q) := \sqrt{\sum_{i=1}^{m} arccos^{2}\lambda_{i}} \nonumber
\end{eqnarray}
$d_{c}$ is a local - but not global - distance function; for instance, if $P$ and $Q$ 
define the same $m$-plane, but with opposite orientations, then $d_{c}(P,Q) = 0$.
The metric $g_{c}$ is generated by $d_{c}$:
\begin{eqnarray}
g_{c}(X,X) := \frac{d}{dt}d_{c}(x(0),x(t))|_{t=0} \nonumber
\end{eqnarray}
where $x=x(t)$ is any smooth path in $\widetilde{G}_{m,n}$ with $\dot{x}(0)= X$.

For our purposes it will be convenient to consider another metric $g_{s}$ on 
$\widetilde{G}_{m,n}$, which is related to spherical geometry. Put 
\begin{eqnarray}
d_{s}(P,Q) := arccos\left( \prod_{i=1}^{m}\lambda_{i} \right) \nonumber
\end{eqnarray}
and let $g_{s}$ be generated by $d_{s}$. It shall be shown below that
$g_{s}$ is, in fact, a well-defined Riemannian metric.
We will prove that $d_{s} \leq d_{c}$. 

\begin{lemma}\label{compare} $g_{s} \leq g_{c}$
\end{lemma} 
 
Next, we consider a natural realization of the metric $g_{s}$. 
The vector space $\Lambda^{m}R^{k}$, where $m\leq k$, possesses a canonical inner-product 
induced from the standard inner-product on $R^{k}$:
\begin{eqnarray} <\xi_{1}\wedge \cdots \wedge \xi_{m}, 
\hspace{0.02in} \zeta_{1}\wedge \cdots \wedge \zeta_{m} >_{\Lambda^{m}R^{k}}
\hspace{0.05in}  & := & \hspace{0.05in} det <\xi_{i},\hspace{0.02in} \zeta_{j}>_{R^{k}} 
\nonumber
\end{eqnarray} 
The unit sphere in $\Lambda^{m}R^{k}$ shall be denoted $S^{\mu}$. 
Consider the submanifold $H$ of $S^{\mu}$
consisting of all elements of the form
$\xi_{1}\wedge \cdots \wedge \xi_{m} \in \Lambda^{m} R^{k}$.
There is a canonical diffeomorphism $\rho: H 
\longrightarrow \widetilde{G}_{m,n}$ defined by
\[ \rho(\xi_{1}\wedge \cdots \wedge \xi_{m}) := (span 
 \{\xi_{1},...,\xi_{m} \}, [ \xi_{1}\wedge \cdots \wedge \xi_{m}] ) \]
 where $[ \xi_{1}\wedge \cdots \wedge \xi_{m}]$ denotes the orientation class of
 $\xi_{1}\wedge \cdots \wedge \xi_{m}$.
 
Let $d_{H}$ be the restriction of the distance function $d_{S^{\mu}}$ 
on the sphere $S^{\mu}$ to 
\[ {\cal D} := \{ (u,v) \in H\times H \hspace*{0.05in} : \hspace*{0.05in}
 <u,v>_{\Lambda^{m} R^{k}} \hspace*{0.02in}
 \mbox{is non-negative} \} \] 
For $(p,q) \in {\cal D}$,
\[ d_{H}(p,q)= d_{S^{\mu}}(p,q) = arccos <p,q>_{\Lambda^{m} R^{k}} =  
 d_{s}(\rho(p),\rho(q))     \]    
Consequently, the  metrics $h$ and $g_{s}$ generated by $d_{H}$ and $d_{s}$, 
respectively, are related by a pull-back: $h = \rho^{*}g_{s}$. 

\begin{lemma}\label{induced}
$h$ is the induced metric on $H$, regarded as a submanifold of the Euclidean space 
$\Lambda^{m}R^{k}$.
\end{lemma}

\begin{corollary}\label{isometry}
$h$ and $g_{s}$ are well-defined and positive definite.
Moreover, $\rho$ is an isometry of Riemannian manifolds $(H, h)$ and
$(\widetilde{G}_{m,n}, g_{s})$. 
\end{corollary}

Let $M$ be an $m$-dimensional, oriented submanifold  of $R^{k}$ and 
let $Z_{1},...,Z_{m}$ be a local, oriented, orthonormal frame 
for $M$.
Define the map $\phi: M \longrightarrow H \subseteq S^{\mu}$ 
by
\[\phi := Z_{1}\wedge \cdots \wedge Z_{m} \]
$\rho\circ\phi:M \longrightarrow \widetilde{G}_{m,n}$ is the (first) Gauss map.

The differential of $\phi$ is 
\begin{eqnarray}
\phi_{*}(X)      & = & \sum_{i=1}^{m} Z_{1}\wedge \cdots \wedge Z_{i-1}\wedge 
                  B(X,Z_{i}) \wedge Z_{i+1} \wedge \cdots \wedge Z_{m} \nonumber
\end{eqnarray} 
Define the homomorphism of vector bundles $B:TM \longrightarrow Hom(TM,TM^{\bot})$ 
by $B(X)(Y) := B(X,Y)$.
We may extend the domain of definition of $B$ by requiring it to act as a derivation on 
tensor products. Then the above equation may be expressed succinctly as 
\begin{eqnarray} \label{differential}
\phi_{*}(X) & = & B(X) \phi 
\end{eqnarray}
In what follows, $Ric$ designates the Ricci curvature of $M$.
\begin{lemma}
\begin{eqnarray}
 \phi^{*}(h)(X,Y) & = & <B(X,Y),trB> - Ric(X,Y) \nonumber
\end{eqnarray}
\end{lemma}

This leads to a characterization of the minimal submanifolds of $R^{k}$.
\begin{corollary}\label{minsub} $M$ is a minimal submanifold of $R^{k}$ if and only if 
\[ \phi^{*}(h) =  - Ric \]
\end{corollary}

\begin{theorem} \label{Euclidean}
Let $M$ be a complete, oriented, minimal submanifold of $R^{k}$ and suppose the 
Gauss map is a bounded embedding. Then $Ric \hspace{0.02in} < \hspace{0.02in}  0$ 
 and
\begin{eqnarray}
\sup Ric & = & 0 \nonumber
\end{eqnarray}
\end{theorem}

\section{The Second Gauss Map}

We suppose that $M$ is an oriented submanifold of $N$, which is an oriented submanifold 
of $R^{k}$.  $B_{M\subseteq N}$ (resp. $B_{N\subseteq E}$) shall denote the 
second fundamental form of $M$ (resp. $N$) viewed as a submanifold of $N$ 
(resp. $R^{k}$). The discussion below proceeds 
in a manner similar to the previous section and so will not be as detailed.

Let $Z_{1},...,Z_{m},V_{1},...,V_{n}$ be a local,  oriented, orthonormal frame for
$N$, where $Z_{1},...,Z_{m}$ is a local oriented frame for $M$, and put $r:=k-n$. 
Define the map $\psi:M\rightarrow S^{\nu}$ by
\[ \psi = V_{1}\wedge\cdots \wedge V_{n} \]
where $S^{\nu}$ is the unit sphere in $\Lambda^{n}R^{k}$.

$H$ shall be the subset of $S^{\nu}$ consisting of elements of the form 
$\xi_{1}\wedge\cdots\wedge \xi_{n}$ and $h$ shall denote the metric on $H$ induced from 
$\Lambda^{n}R^{k}$.  The diffeomorphism $\rho: H 
\longrightarrow \widetilde{G}_{n,r}$ defined by
\[ \rho(\xi_{1}\wedge \cdots \wedge \xi_{n}) := (span 
 \{\xi_{1},...,\xi_{n} \}, [ \xi_{1}\wedge \cdots \wedge \xi_{n}] ) \]
 is an isometry of Riemannian manifolds $(H, h)$ and
$(\widetilde{G}_{n,r}, g_{s})$. $\rho\circ\psi:M\longrightarrow \widetilde{G}_{n,r}$
is the {\it second Gauss map}.
                   
\begin{lemma}\label{seclemma} 
\begin{eqnarray}
\psi^{*}(h)(X,Y) 
& = & <B_{M\subseteq N}(X,Y), trB_{M\subseteq N}> -Ric_{M}(X,Y) + \nonumber \\
&   & \sum_{i=1}^{n} <B_{N\subseteq E}(X,Y), B_{N\subseteq E}(V_{i},V_{i})> + \nonumber \\
&   & \sum_{i=1}^{m} <R_{N}(X,Z_{i})Y,Z_{i}> - \sum_{i=1}^{n} <R_{N}(X,V_{i})Y,V_{i}>   
 \nonumber     
\end{eqnarray} 
\end{lemma}

\begin{corollary}
If $M^{m}$ is a minimal submanifold of $S^{k-1}$, the unit hypersphere of $R^{k}$ then
\begin{eqnarray}
\psi^{*}(h) 
& = & (m-1)<,>_{M} -Ric \nonumber
\end{eqnarray} 
\end{corollary}
The proof of the theorem below is essentially contained in the proof of  
Theorem \ref{Euclidean}. 

\begin{theorem} Let $M^{m}$ be a complete, oriented, minimal submanifold of 
a unit sphere
and suppose the second Gauss map is a bounded embedding. 
Then $Ric\hspace*{0.02in} < \hspace*{0.02in} m-1$. 
Furthermore, $M$ is compact  if and only if 
\begin{eqnarray}
\sup Ric & < & m-1 \nonumber
\end{eqnarray} 
\end{theorem}

\end{document}